\documentclass[11pt]{amsart}

\usepackage[margin=1.2in]{geometry}
\usepackage{amsmath,amssymb}
\usepackage{array}
\usepackage{booktabs}
\usepackage{graphicx}
\usepackage{newtxtext,newtxmath}
\usepackage{hyperref}
\usepackage{xcolor}
\usepackage{tikz}
\usetikzlibrary{arrows.meta,calc,decorations.markings,positioning}

\hypersetup{
  colorlinks=true,
  linkcolor=blue!55!black,
  citecolor=blue!55!black,
  urlcolor=blue!55!black
}

\newcommand{\bd}{\partial}

\definecolor{tirelight}{RGB}{225,228,231}
\definecolor{tiregray}{RGB}{91,96,103}
\definecolor{rimred}{RGB}{185,42,45}
\definecolor{loopblue}{RGB}{30,103,166}
\definecolor{facegold}{RGB}{244,208,118}
\definecolor{faceteal}{RGB}{137,206,194}

\newcommand{\torusbody}{%
  \draw[draw=tiregray, double=tirelight, double distance=12mm, line width=.8pt,
    line cap=butt]
    (0:1.45) arc (0:180:1.45);

  \draw[draw=tiregray, double=tirelight!86!gray, double distance=12mm,
    line width=.8pt, line cap=butt]
    (180:1.45) arc (180:360:1.45);

  \foreach \torusx in {-1.45,1.45}{%
    \begin{scope}[shift={(\torusx,0)}]
      \pgftransformresetnontranslations
      \fill[tirelight!86!gray]
        (-6.25mm,-.3mm) -- (-6.25mm,0)
        arc[start angle=180,end angle=0,radius=6.25mm]
        -- (6.25mm,-.3mm) -- cycle;
      \draw[tiregray, line width=.8pt] (6.15mm,0)
        arc[start angle=0,end angle=180,radius=6.15mm];
      \draw[tiregray, line width=.8pt, dashed] (-6.15mm,0)
        arc[start angle=180,end angle=360,radius=6.15mm];
    \end{scope}
  }%
}

\title{How Many Sides Does a Tire Have?}
\author[Vanja Stojanovi\'c]{Vanja Stojanovi\'c\\
{\small University of Ljubljana, Faculty of Mathematics and Physics}\\
{\small Jadranska ulica 21}}
\address{University of Ljubljana, Faculty of Mathematics and Physics,
Jadranska ulica 21}
\date{September 2026}
\subjclass[2020]{00A05, 55N35, 57K20}
\keywords{Topology, torus, solid torus, cellular chains, boundary operator}

\begin{document}

\begin{abstract}
How many sides does a tire have: two, three, or only one?  The question is a
joke only until we ask what the word ``side'' is supposed to measure.  This
article follows three reasonable answers into one small piece of topology: a
loop on a tire can have no endpoints and still be a border of a chosen region.
The point is not to find the official number of sides, but to see why the count
changes when the model changes.
\end{abstract}

\maketitle

\section{A question with too many good answers}

Ask a few people how many sides a tire has.  Someone will say two and point to
the sidewalls.  Someone else will count the sidewalls and the tread and say
three.  A third person may run a finger over the rubber, notice that there is no
edge to fall off, and insist that the whole surface is one piece.  Each answer
is responding to a different version of the question.

This article is about that change of version.  We will not ask for the correct
number of sides.  We will ask what must be specified before such a number has a
meaning.  Is a side a part with a familiar name, a connected patch of rubber, or
a face in a chosen decomposition?  Once that choice is made, the counting can
begin.

That shift from object to definition is a basic topological move.  The tire
alone does not settle the count; the structure we place on it matters too.  The
payoff will be a distinction that is easy to blur in ordinary language: a curve
can have no boundary of its own while serving as the boundary of a region.

\section{The tire in three-dimensional space}

\begin{figure}[ht]
\centering
\includegraphics[width=.42\linewidth]{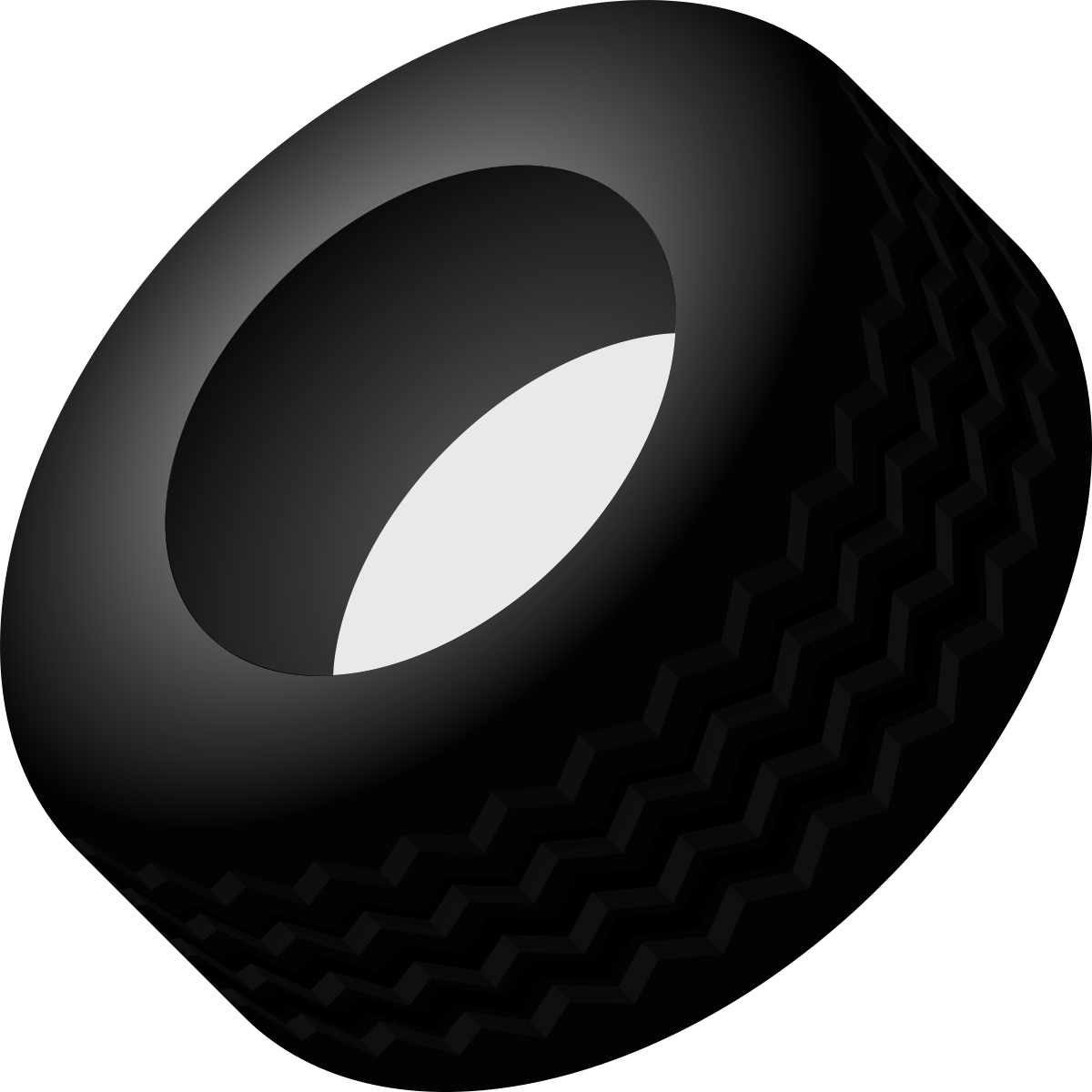}
\caption{A car tire is a thick three-dimensional object.  Topology ignores the
tread pattern, branding, and small geometric details while keeping the circular
hole and the solid rubber body.  Public-domain illustration by barretr,
\href{https://openclipart.org/detail/188729/tire}{Openclipart}.}
\label{fig:car-tire}
\end{figure}

A car tire is a three-dimensional body, not merely a doughnut-shaped curve.
Its tread, sidewalls, and deep central opening are visible in
Figure~\ref{fig:car-tire}.

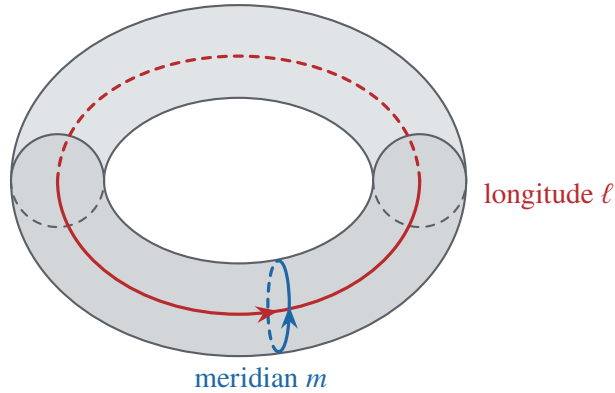
\begin{figure}[ht]
\centering
\begin{tikzpicture}[line cap=round, line join=round]
  \begin{scope}[xscale=1.65,yscale=1.18]
    \torusbody

    \draw[loopblue, very thick, dashed]
      (.32,-1.92) to[out=180,in=180,looseness=.29] (.32,-.9);

    \draw[rimred, very thick, dashed]
      (1.45,0) to[out=92,in=88,looseness=1.65] (-1.45,0);
    \draw[rimred, very thick,
      postaction={decorate}, decoration={markings,
      mark=at position .58 with {\arrow{Stealth}}}]
      (-1.45,0) arc[start angle=180,end angle=360,
        x radius=1.45,y radius=1.5];

    \draw[loopblue, very thick,
      postaction={decorate}, decoration={markings,
      mark=at position .50 with {\arrow{Stealth}}}]
      (.32,-1.92) to[out=0,in=0,looseness=.29] (.32,-.9);
  \end{scope}

  \node[rimred, inner sep=1.5pt] at (4.1,-.20)
    {longitude \(\ell\)};
  \node[loopblue, inner sep=1.5pt] at (.3,-2.6)
    {meridian \(m\)};
\end{tikzpicture}
\caption{The boundary torus in a standard two-dimensional perspective drawing.
The meridian \(m\) circles the tube, while the longitude \(\ell\) follows the
tire's rolling direction.  Dashed arcs pass behind the visible surface.}
\label{fig:torus-coordinates}
\end{figure}

Ignoring the tread pattern and the slight flattening at the road, the standard
topological model is the \emph{solid torus}
\[
  V=S^1\times D^2.
\]
Here \(S^1\) means a circle and \(D^2\) means a filled disk.  The product sign
\(\times\) asks us to attach one copy of the disk to every point of the circle.
More physically, imagine carrying the disk once around the circle: the disk
sweeps out all the rubber.  The surface of the rubber is
\[
  T=\bd V=S^1\times S^1,
\]
the ordinary torus shown in Figure~\ref{fig:torus-coordinates}.

The two colored loops are worth naming.  A \emph{meridian} goes around the
cross-section of the tube.  A \emph{longitude} goes around the central hole.
These are the two independent circular directions in \(S^1\times S^1\).

As a topological space, \(T\) is connected and has no boundary:
\[
  \bd T=\varnothing.
\]
The symbol \(\bd\) here means the ordinary geometric boundary: \(T=\bd V\) is
the surface enclosing the solid tire, while \(\bd T=\varnothing\) says that this
surface has no edge of its own.  This does not forbid us from drawing edges on
it.  A map painted on a globe has country borders even though the sphere itself
has no boundary.  We will shortly use \(\bd_k\) for a related but different
idea: the boundary of one of the pieces in such a map.

\section{What a boundary operator records}

Before introducing notation, here is the issue in plain language.  The whole
surface of the tire has no edge where the rubber surface stops.  But if we draw
or choose a region on that surface, the region may have a border.  The notation
below is a bookkeeping device for keeping those two uses of ``boundary''
separate.

To turn painted regions into mathematics, we use \emph{cells}.  A 0-cell is a
vertex, a 1-cell is an edge, and a 2-cell is a disk-shaped face.  For each
dimension \(k\), the notation \(C_k\) means all finite formal sums of the
\(k\)-cells.  ``Formal sum'' is only bookkeeping.  For example, an expression
such as \(a+b-c\) says that the edges \(a\) and \(b\) are being counted in their
chosen directions, while \(c\) is being counted in the opposite direction.

Choosing an \emph{orientation} means choosing a direction along each edge and
around each face.  The cellular boundary operator
\[
  \bd_k:C_k\longrightarrow C_{k-1}
\]
sends each oriented piece to the signed collection of pieces surrounding it.
It also respects addition: the boundary of a sum is the sum of the boundaries.
For an oriented line segment,
\[
  \bd_1(\text{segment})=(\text{end})-(\text{start}).
\]
A closed loop has no unmatched start or end, so its 1-dimensional boundary is
zero.  A chain whose boundary is zero is called a \emph{cycle}.

Let \(e\) be a marked longitude representing an inner-rim line on our tire.  In
both models below,
\[
  \bd_1e=0.
\]
That equation says only that \(e\) is a closed cycle.  It does \emph{not} tell us
whether \(e\) is used as an edge of a face or region.  For that we must look one
dimension higher, at \(\bd_2\).

\section{Model A: the rim belongs to two regions}

One more familiar model will help us see the seams.  Start with a square and
glue its top edge to its bottom edge, matching the arrows; this first makes a
cylinder.  Then glue the two remaining edges to each other.  The result is a
torus.  Thus a path that leaves one side of the square immediately reappears at
the matching point on the opposite side.

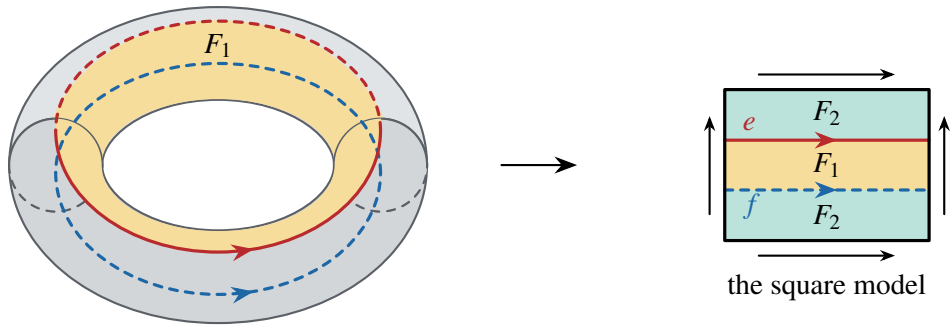
\begin{figure}[ht]
\centering
\begin{tikzpicture}[line cap=round, line join=round]
  \begin{scope}[xscale=1.48,yscale=1.02,shift={(-1.55,0)}]
    \torusbody

    \fill[facegold!72]
      (-1.56,0)
      arc[start angle=180,end angle=0,x radius=1.56,y radius=1.88]
      (1.04,0)
      arc[start angle=0,end angle=180,x radius=1.04,y radius=.857]
      -- cycle;

    \foreach \torusx in {-1.45,1.45}{%
      \begin{scope}[shift={(\torusx,0)}]
      \pgftransformresetnontranslations
      \fill[tirelight!86!gray] (0,0) circle[radius=6.15mm];
      \draw[tiregray, line width=.8pt] (6.15mm,0)
        arc[start angle=0,end angle=180,radius=6.15mm];
      \draw[tiregray, line width=.8pt, dashed] (-6.15mm,0)
        arc[start angle=180,end angle=360,radius=6.15mm];
    \end{scope}
    }%

    \fill[facegold!72]
      (-1.04,0)
      arc[start angle=180,end angle=360,x radius=1.04,y radius=.857]
      arc[start angle=180,end angle=90,x radius=.415,y radius=.603]
      -- (1.45,.45)
      arc[start angle=0,end angle=-180,x radius=1.45,y radius=1.58]
      -- (-1.455,.603)
      arc[start angle=90,end angle=0,x radius=.415,y radius=.603]
      -- cycle;

    \begin{scope}[shift={(0,-.35)}]
      \draw[loopblue, very thick, dashed]
        (1.45,0.25) arc[start angle=0,end angle=180,
          x radius=1.45,y radius=1.42];
      \draw[loopblue, very thick, dashed,
        postaction={decorate}, decoration={markings,
        mark=at position .58 with {\arrow{Stealth}}}]
        (-1.45,0.25) arc[start angle=180,end angle=360,
          x radius=1.45,y radius=1.58];
    \end{scope}
    \begin{scope}[shift={(0,.45)}]
      \draw[rimred, very thick, dashed]
        (1.45,0) arc[start angle=0,end angle=180,
          x radius=1.45,y radius=1.42];
      \draw[rimred, very thick,
        postaction={decorate}, decoration={markings,
        mark=at position .58 with {\arrow{Stealth}}}]
        (-1.45,0) arc[start angle=180,end angle=360,
          x radius=1.45,y radius=1.58];
    \end{scope}
    \node at (0,1.55) {\(F_1\)};
  \end{scope}

  \draw[-{Stealth[length=3mm]}, thick] (1.45,0) -- (2.45,0);

  \begin{scope}[shift={(5.75,0)}]
    \fill[faceteal!58] (-1.35,-1.0) rectangle (1.35,1.0);
    \fill[facegold!72] (-1.35,-.33) rectangle (1.35,.33);
    \draw[rimred, very thick,
      postaction={decorate}, decoration={markings,
      mark=at position .55 with {\arrow{Stealth}}}]
      (-1.35,.33) -- (1.35,.33);
    \draw[loopblue, very thick, dashed,
      postaction={decorate}, decoration={markings,
      mark=at position .55 with {\arrow{Stealth}}}]
      (-1.35,-.33) -- (1.35,-.33);
    \draw[very thick] (-1.35,-1.0) rectangle (1.35,1.0);
    \draw[-{Stealth[length=2mm]}, thick] (-.9,1.2) -- (.9,1.2);
    \draw[-{Stealth[length=2mm]}, thick] (-.9,-1.2) -- (.9,-1.2);
    \draw[-{Stealth[length=2mm]}, thick] (-1.55,-.65) -- (-1.55,.65);
    \draw[-{Stealth[length=2mm]}, thick] (1.55,-.65) -- (1.55,.65);
    \node at (0,0) {\(F_1\)};
    \node at (0,.69) {\(F_2\)};
    \node at (0,-.69) {\(F_2\)};
    \node[rimred, anchor=west] at (-1.25,.55) {\(e\)};
    \node[loopblue, anchor=west] at (-1.25,-.55) {\(f\)};
    \node[align=center] at (0,-1.62) {the square model};
  \end{scope}
\end{tikzpicture}
\caption{Model A uses two parallel seams.  In the square model, the top and
bottom teal bands join to form the annular chain \(F_2\).  The chosen inner rim
\(e\) is genuinely a shared edge of the two named regions.}
\label{fig:model-a}
\end{figure}

Suppose we paint the torus in two colors and decree that the color changes along
the inner rim \(e\).  There is a subtlety: one longitude by itself does not cut
a torus into two pieces.  We need a second parallel seam \(f\).  Together, the
two seams divide the surface into two annuli, or cylinder-shaped regions, as
suggested in Figure~\ref{fig:model-a}.  An annulus is not itself a disk-shaped
cell, but we can split it into disk-shaped faces by drawing one additional
crosswise edge.  Below, \(F_1\) and \(F_2\) denote the chains obtained by adding
the faces in the two annular regions.  The new crosswise edges occur twice with
opposite signs and therefore cancel from each region's total boundary.

With compatible orientations, the relevant part of the cellular boundary is
\[
  \bd_2F_1=e-f,
  \qquad
  \bd_2F_2=f-e.
\]
The signs reverse because the pieces on the two sides traverse each shared edge
in opposite directions.  In particular, \(e\) occurs with a nonzero coefficient in both
boundaries.  The boundaries cancel when the two region-chains are added:
\[
  \bd_2(F_1+F_2)=(e-f)+(f-e)=0,
\]
just as they must when the two regions are reassembled into the boundaryless
torus.

In this model, calling \(e\) a border is completely reasonable.  Crossing it
moves us from one named region to the other.  The second seam is not an annoying
technicality, but a feature that makes the two-region picture possible on a
torus.

\section{Model B: the rim is only a loop}

\begin{figure}[ht]
\centering
\begin{tikzpicture}[line cap=round, line join=round]
  \begin{scope}[xscale=1.48,yscale=1.06,shift={(-1.55,0)}]
    \torusbody

    \draw[loopblue, very thick, dashed]
      (.32,-1.97) to[out=180,in=180,looseness=.3] (.32,-0.85);

    \draw[rimred, very thick, dashed]
      (1.45,0) to[out=92,in=88,looseness=1.65] (-1.45,0);
    \draw[rimred, very thick,
      postaction={decorate}, decoration={markings,
      mark=at position .58 with {\arrow{Stealth}}}]
      (-1.45,0) arc[start angle=180,end angle=360,
        x radius=1.45,y radius=1.5];

    \draw[loopblue, very thick,
      postaction={decorate}, decoration={markings,
      mark=at position .50 with {\arrow{Stealth}}}]
      (.32,-1.97) to[out=0,in=0,looseness=.3] (.32,-.85);
  \end{scope}
  \draw[-{Stealth[length=3mm]}, thick] (1.35,0) -- (2.1,0);

  \begin{scope}[shift={(5.15,0)}]
 \fill[tirelight] (-1.35,-1.0) rectangle (1.35,1.0);
    \draw[rimred, very thick,
      postaction={decorate}, decoration={markings,
      mark=at position .55 with {\arrow{Stealth}}}]
      (-1.35,-.23) -- (1.35,-.23);
    \draw[loopblue, very thick,
      postaction={decorate}, decoration={markings,
      mark=at position .55 with {\arrow{Stealth}}}]
      (.25,-1.0) -- (.25,1.0);
   
    \draw[very thick] (-1.35,-1.0) rectangle (1.35,1.0);
    \draw[-{Stealth[length=2mm]}, thick] (-.9,1.2) -- (.9,1.2);
    \draw[-{Stealth[length=2mm]}, thick] (-.9,-1.2) -- (.9,-1.2);
    \draw[-{Stealth[length=2mm]}, thick] (-1.55,-.65) -- (-1.55,.65);
    \draw[-{Stealth[length=2mm]}, thick] (1.55,-.65) -- (1.55,.65);
    \node[rimred, anchor=south west] at (-1.25,-.18) {\(e\)};
    \node[loopblue, anchor=south west] at (.30,.5) {\(m\)};
    \node[align=center] at (0,-1.63)
      {opposite sides are identified};
  \end{scope}
\end{tikzpicture}
\caption{Model B on the embedded torus and on its square model.  The blue
meridian crosses the red longitude once.  Because opposite sides of the square
are identified, both paths close up.}
\label{fig:model-b}
\end{figure}
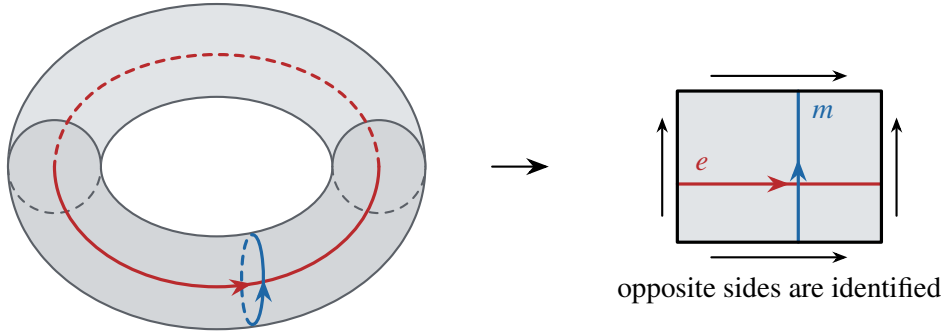

Now erase the second seam.  The remaining red longitude \(e\) is still closed,
so \(\bd_1e=0\), but it no longer separates two regions.  Figure~\ref{fig:model-b}
shows a quick way to see this.  The blue meridian crosses \(e\) once and returns
to its starting point.  If the red curve really separated the surface into two
regions, any closed trip would have to cross it an even number of times: once
to leave a region and once to come back.

Cutting a torus along one longitude makes a cylinder, and a cylinder is still
connected.  This is the formal content of saying that \(e\) is
\emph{nonseparating}.  Thus the red loop may be visually prominent without
being a border between two global pieces.

The familiar analogy with the equator is slightly dangerous here.  An equator
\emph{does} separate a sphere into northern and southern hemispheres, whether
or not we color them.  A longitude on a torus behaves differently: removing it
leaves the rest of the torus in one piece.  The hole in the torus is what allows
the blue route in Figure~\ref{fig:model-b} to cross the red loop only once.

\section{So which boundary is zero?}

The phrase ``the boundary of the rim'' can now be unpacked.  The comparison is
not
\[
  \bd e\ne0 \quad\hbox{versus}\quad \bd e=0.
\]
The rim is a closed loop, hence a cycle, in both models, so \(\bd_1e=0\) in
both.  What changes is whether \(e\) appears in the boundary of a
2-dimensional region:

\begin{center}
\begin{tabular}{>{\raggedright\arraybackslash}p{.16\linewidth}
                >{\raggedright\arraybackslash}p{.47\linewidth}
                >{\raggedright\arraybackslash}p{.24\linewidth}}
\toprule
Model & What the boundary maps say & What \(e\) does \\
\midrule
A & \(\bd_1e=0\), while \(e\) occurs in \(\bd_2F_1\) and \(\bd_2F_2\)
  & shared cellular border \\
B & \(\bd_1e=0\), with no two regions having \(e\) as their interface
  & nonseparating marked loop \\
\bottomrule
\end{tabular}
\end{center}

This distinction is one reason algebraic topology keeps track of dimension so
carefully.  The map \(\bd_1\) asks for the endpoints of a curve.  The map
\(\bd_2\) asks for the edges of a face or region.  A loop can have no endpoints
and still occur as part of a region's boundary.  There is no contradiction
because the two boundary maps answer different questions.  This is the
mathematical point hidden in the original tire question.

\section{How many sides, then?}

There is no single number until we specify a convention.  If ``sides'' means
named functional regions of a manufactured tire, two sidewalls plus a tread may
give three.  If it means the named regions in Model A, our chosen decomposition
gives two annular bands, each bounded by the same pair of seams.  If it means
connected components of the unmarked rubber surface, there is just one.

One should also resist saying simply that a torus is ``one-sided.''  In the
ordinary tire-shaped torus sitting in three-dimensional space, an arrow
perpendicular to the surface can point inward or outward.  In this geometric
sense the torus is two-sided, unlike a M\"obius strip, whose twist exchanges the
two local choices.  That is a different use of \emph{side} again.

So the best answer is a question: \emph{which sides are we counting?}  What
began as wordplay has exposed a central habit of topology.  Before counting
pieces, first decide what makes two pieces different.  Before speaking of a
boundary, say whose boundary and in which dimension.  The tire has not changed;
our chosen model has.

\renewcommand{\refname}{Further Reading}

\end{document}